\title{Model Validation and Vibration Analysis of Cosserat Plates}

\author{
Lev Steinberg$\mathbf{^{1}}$ \\
$^\mathbf{1}$ Department of Mathematical Sciences \\
University of Puerto Rico at Mayag\"{u}ez,\\
Mayag\"{u}ez, Puerto Rico 00681-9018, USA\\
Roman Kvasov$\mathbf{^{2}}$\\
$^\mathbf{2}$ Department of Mathematics\\
University of Puerto Rico at Aguadilla \\
Aguadilla, Puerto Rico 00604, USA
}
\date{\today}

\documentclass[12pt]{article}
\usepackage{geometry}
\usepackage{setspace}
\usepackage{graphicx}
\usepackage{amssymb}
\usepackage{amsfonts}
\usepackage{amsmath}
\usepackage{array}
\usepackage{booktabs,tabularx}
\usepackage{float}
\usepackage{multicol}
\usepackage[utf8]{inputenc}
\usepackage[T1]{fontenc}

\begin{document}

\maketitle

\begin{abstract}
In this paper we present the validation of our recently published mathematical model for the dynamics of Cosserat elastic plates. The validation is based on the comparison with the exact solution of the 3-dimensional Cosserat elastodynamics. The preliminary computations of eigenfrequencies show the high agreement with the exact values. The computations allow us to detect the splitting of the frequencies of vibrations (micro vibration) depending on the orientation of micro elements. This provided us with a powerful tool for distinguishing between the frequencies of the micro and macro vibrations of the plate. 
\end{abstract}

\section{Introduction}

The theory of asymmetric elasticity introduced in 1909 by the Cosserat brothers \cite{Cosserat1909} gave rise to a variety of Cosserat plate theories. In 1960s Green and Naghdi specialized their general theory of Cosserat surface to obtain the linear Cosserat plate \cite{Green1966}, while independently Eringen proposed a complete theory of plates in the framework of Cosserat elasticity \cite{Eringen1967}.

The first theory of Cosserat plates based on the Reissner plate theory was developed in \cite{Steinberg2010} and its finite element modeling is provided in \cite{Kvasov2011}. The parametric theory of Cosserat plate, presented by the authors in \cite{Steinberg2013}, includes some additional assumptions leading to the introduction of the splitting parameter. This guaranteed the highest level of approximation to the original three-dimentional problem. The parametric theory produces the equilibrium equations, constitutive relations, and the optimal value of the minimization of the elastic energy of the Cosserat plate. The paper \cite{Steinberg2013} also provides the analytical solutions of the presented plate theory and the three-dimensional Cosserat elasticity for simply supported rectangular plate. The comparison of these solutions showed that the precision of the developed Cosserat plate theory is similar to the precision of the classical plate theory developed by Reissner \cite{Reissner1944}, \cite{Reissner1945}.

The numerical modeling of bending of simply supported rectangular plates is given in \cite{Kvasov2013}. We developed the Cosserat plate field equations and a rigorous formula for the optimal value of the splitting parameter. The solution of the Cosserat plate was shown to converge to the Reissner plate as the elastic asymmetric parameters tend to zero. The Cosserat plate theory demonstrates the agreement with the size effect, confirming that the plates of smaller thickness are more rigid than is expected from the Reissner model. The modeling of Cosserat plates with simply supported rectangular holes is also provided.

The extension of the static model of Cosserat elastic plates to the dynamic problems is presented in \cite{Steinberg2015}. The computations predict a new kind of natural frequencies associated with the material microstructure and were shown to be compatible with the size effect principle reported in \cite{Kvasov2013} for the Cosserat plate bending.

The numerical study of Cosserat elastic plate deformation based on the parametric theory of Cosserat plates using the Finite Element Method is presented in \cite{Kvasov2017}. The paper discusses the existence and uniqueness of the weak solution, convergence of the proposed FEM and its numerical validation by estimating the order of convergence. The Finite Element analysis of clamped Cosserat plates of different shapes under different loads is also provided. The numerical analysis of plates with circular holes shows that the stress concentration factor around the hole is less than the classical value, and smaller holes exhibit less stress concentration as would be expected on the basis of the classical elasticity.

The current article represents an extension of the paper \cite{Steinberg2015} for different shapes and orientations of micro-elements incorporated into the Cosserat plates. It is based on the generalized variational principle for elastodynamics and includes a non-diagonal rotatory inertia tensor. The numerical computations of the plate free vibration showed the existence of some additional high frequencies of micro-vibrations depending on the orientation of micro-elements. The comparison with three-dimensional Cosserat elastodynamics shows a high agreement with the exact values of the eigenvalue frequencies.

\section{Cosserat Linear Elastodynamics}

\subsection{Fundamental Equations}

The Cosserat linear elasticity balance laws are
\begin{eqnarray}
\mathbf{div}\boldsymbol{\sigma } &=&\frac{\partial \mathbf{p}}{\partial t},  \label{equilibrium_equations 0} \\
\boldsymbol{\varepsilon}\cdot \boldsymbol{\sigma}+\mathbf{div} \boldsymbol{\mu} &=&\frac{\partial \mathbf{q}}{\partial t}, \label{equilibrium_equations}
\end{eqnarray}
where the $\boldsymbol{\sigma}$ is the stress tensor, $\boldsymbol{\mu}$ the couple stress tensor, $\mathbf{p}=\rho \frac{\partial \mathbf{u}}{\partial t}$ and $\mathbf{q}=\mathbf{J}{\frac{\partial \boldsymbol{\phi}}{\partial t}}$ are the linear and angular momenta, $\rho$ and $\mathbf{J}$ are the material density and the rotatory inertia characteristics, $\boldsymbol{\varepsilon}$ is the Levi-Civita tensor.

We will also consider the constitutive equations in the following form \cite{Nowacki1986}: 
\begin{eqnarray}
\mathbf{\boldsymbol{\sigma}} &=& (\mu +\alpha ) \boldsymbol{\gamma }+(\mu-\alpha) \mathbf{\boldsymbol{\gamma }}^{T}+\lambda \mathbf{(tr%
\boldsymbol{\gamma})1},  \label{Hooke's_law 1} \\
\mathbf{\boldsymbol{\mu }} &=& (\gamma +\epsilon ) \boldsymbol{\chi }+(\gamma -\epsilon) \boldsymbol{\chi }^{T}+\beta \mathbf{(tr\boldsymbol{\chi}%
)1},  \label{Hooke's_law 1A}
\end{eqnarray}
and the kinematic relations in the form 
\begin{equation}
\mathbf{\boldsymbol{\gamma}=}\left( \mathbf{\nabla u}\right) ^{T}+\boldsymbol{\varepsilon}\cdot \boldsymbol{\phi }\text{ and }\boldsymbol{\chi
}=\left( \mathbf{\nabla \boldsymbol{\phi}}\right) ^{T},
\label{kinematic formulas}
\end{equation}
Here $\mathbf{u}$ and $\boldsymbol{\phi}$ represent the displacement and rotation vectors, $\boldsymbol{\gamma}$ and $\boldsymbol{\chi}$ represent the strain and torsion tensors, $\mu$, $\lambda $ are the Lam\'{e} parameters and $\alpha$, $\beta$, $\gamma$, $\epsilon$ are the Cosserat elasticity parameters.

The constitutive equations (\ref{Hooke's_law 1}) - (\ref{Hooke's_law 1A}) can be written in the reverse form \cite{Steinberg2010}:
\begin{eqnarray}
\boldsymbol{\gamma } &=& (\mu ^{\prime }+\alpha ^{\prime })%
\boldsymbol{\sigma }+(\mu ^{\prime }-\alpha ^{\prime })\boldsymbol{\sigma }%
^{T}+\lambda ^{\prime }\mathbf{(tr\boldsymbol{\sigma})1},
\label{Hooke's_law 2} \\
\boldsymbol{\chi } &=& (\gamma ^{\prime }+\epsilon ^{\prime })%
\boldsymbol{\mu }+(\gamma ^{\prime }-\epsilon ^{\prime })\boldsymbol{\mu }%
^{T}+\beta^{\prime } \mathbf{(tr\boldsymbol{\mu})1}.  \label{Hooke's_law 2A}
\end{eqnarray}
where $\mu^{\prime}=\frac{1}{4\mu}$, $\alpha ^{\prime}=\frac{1}{4\alpha}$, $\gamma^{\prime}=\frac{1}{4\gamma }$, $\epsilon ^{\prime }=\frac{1}{4\epsilon }$, $\lambda ^{\prime }=\frac{-\lambda }{6\mu (\lambda +\frac{2\mu}{3})}$ and $\beta ^{\prime }=\frac{-\beta }{6\mu (\beta +\frac{2\gamma}{3})}$.

We will consider the boundary conditions provided in  \cite{Steinberg2015}:
\begin{eqnarray}
\mathbf{u}=\mathbf{u}_{\mathbf{0}}, \text{ } \boldsymbol{\phi} =
\boldsymbol{\phi }_{\mathbf{0}}\text{ on }\mathcal{G}_{1}^{t}=\partial
B_{0}\backslash \partial B_{\sigma }\times \left[ t_{0},t\right] ,
\label{Bound conditions} && \\
\boldsymbol{\sigma}_{\mathbf{n}}=\boldsymbol{\sigma} \cdot \mathbf{n}=\boldsymbol{\sigma }_{\mathbf{0}},\text{ on }
\mathcal{G}_{2}^{t}=\partial B_{\sigma }\times \left[ t_{0},t\right], &&\\
\boldsymbol{\mu}_{\mathbf{n}}=\boldsymbol{\mu }\cdot \mathbf{n} =\boldsymbol{\mu}_{\mathbf{0}}\text{ on }
\mathcal{G}_{2}^{t}=\partial B_{\sigma }\times \left[ t_{0},t\right],  &&
\label{Bound conditions a}
\end{eqnarray}%
and initial conditions 
\begin{eqnarray}
\mathbf{u}\left( x,0\right) = \mathbf{U}_{\mathbf{0}}\mathbf{,}%
\text{ }\mathbf{\phi \left( x,0\right) =\Phi }_{\mathbf{0}}\text{ in } B_{0}, &&\\
\mathbf{\dot{u}}\left( x,0\right) = \mathbf{\bar{U}}_{\mathbf{0}}%
\mathbf{,}\text{ }\mathbf{\dot{\phi}\left( x,0\right) =\bar{\Phi}}_{\mathbf{0%
}}\text{ in }B_{0}, &&
\end{eqnarray}%
where $t_{0}$ and $t$ are the initial and terminal time, $\mathbf{u}_{\mathbf{0}}$ and $\boldsymbol{\phi }_{\mathbf{0}}$ are
prescribed on $\mathcal{G}_{1}$, $\boldsymbol{\sigma }_{\mathbf{0}}$ and $\boldsymbol{\mu }_{\mathbf{0}}$\ on $\mathcal{G}_{2},$ and $\mathbf{n}$ is the unit vector normal to the boundary $\partial B_{0}$ of the elastic body $B_{0}$.

\vspace{0.2in}

\subsection{Cosserat Elastic Energy}

\vspace{0.2in}

The strain stored energy $U_{C}$ of the body $B_{0}$ is defined by the integral \cite{Nowacki1986}: 
\begin{equation}
{\large U}_{C}=\int_{B_{0}}\emph{W}\left\{ \mathbf{\boldsymbol{\gamma}, %
\boldsymbol{\chi}}\right\} dv,  \label{free_energy_strain}
\end{equation}
where 
\begin{eqnarray*}
\emph{W}\left\{\boldsymbol{\gamma},\boldsymbol{\chi}\right\} &=&\frac{\mu+\alpha}{2}\gamma_{ij}\gamma_{ij}+\frac{\mu-\alpha}{2}\gamma_{ij}\gamma_{ji} \\ \notag
&&+\frac{\lambda}{2}\gamma_{kk}\gamma_{nn} +\frac{\gamma +\epsilon }{2}\chi _{ij}\chi _{ij} \\ \notag
&&+\frac{\gamma-\epsilon}{2}\chi_{ij}\chi_{ji}+\frac{\beta}{2}\chi _{kk}\chi_{nn}, \notag
\end{eqnarray*}
is non-negative and the relations (\ref{Hooke's_law 1}) - (\ref{Hooke's_law
1A}) can be written in the form \cite{Steinberg2015}: 
\begin{equation}
\boldsymbol{\sigma} = \nabla _{\boldsymbol{\gamma }}\emph{W}\text{ and }%
\boldsymbol{\mu}= \nabla _{\boldsymbol{\chi }}\emph{W}.  \label{Const-d 2}
\end{equation}

The stress energy is given as 
\begin{equation}
{\large U}_{K}=\int_{B_{0}}\Phi \left\{ \boldsymbol{\sigma},\boldsymbol{\mu}%
\right\} dv,
\label{free_energy_stress}
\end{equation}
where 
\begin{eqnarray*}
\Phi \left\{\boldsymbol{\sigma},\boldsymbol{\mu}\right\} &=&\frac{\mu^{\prime}+\alpha ^{\prime }}{2}\sigma _{ij}\sigma _{ij}+\frac{\mu ^{\prime
}-\alpha ^{\prime }}{2}\sigma _{ij}\sigma _{ji} \notag \\
&&+\frac{\lambda ^{\prime }}{2}%
\sigma _{kk}\sigma _{nn} +\frac{\gamma ^{\prime }+\epsilon ^{\prime }}{2}\mu _{ij}\mu _{ij} \notag \\
&&+\frac{\gamma ^{\prime }-\epsilon ^{\prime }}{2}\mu _{ij}\mu _{ji}+\frac{\beta
^{\prime }}{2}\mu _{kk}\mu _{nn},  
\end{eqnarray*}
and the relations (\ref{Hooke's_law 2}) - (\ref{Hooke's_law 2A}) can be written as \cite{Steinberg2015} 
\begin{equation}
\boldsymbol{\gamma} = \frac{\partial \Phi }{\partial \boldsymbol{\sigma}},%
\text{ and }\boldsymbol{\chi} = \frac{\partial \Phi }{\partial %
\boldsymbol{\mu }}.  \label{Const1-1B}
\end{equation}

We consider the work done by the stresses $\boldsymbol{\sigma}$ and $%
\boldsymbol{\mu}$ over the strains $\boldsymbol{\gamma}$ and $%
\boldsymbol{\chi}$ as in \cite{Nowacki1986} 
\begin{equation}
{\large U}=\int_{B_{0}}\left[ \boldsymbol{\sigma}\cdot \boldsymbol{\gamma}+%
\boldsymbol{\mu}\cdot \boldsymbol{\chi}\right] dv  \label{free energy}
\end{equation}%
and 
\begin{equation}
U=U_{K}=U_{C}
\end{equation}

The stored kinetic energy $T_{C}$ is defined as 
\begin{equation}
T_{C}=\int_{B_{0}}\Upsilon _{C}dv=\frac{1}{2} \int_{B_{0}}\left( \rho \left( 
\frac{\partial \mathbf{u}}{\partial t}\right)^{2} + \mathbf{J}\left( {\frac{%
\partial \boldsymbol{\phi }}{\partial t}}\right)^{2}\right) dv,
\end{equation}

The kinetic energy $T_{K}$ is given as 
\begin{equation}
T_{K}=\int_{B_{0}}\Upsilon _{K}\{\mathbf{p},\mathbf{q}\}dv=\frac{1}{2}
\int_{B_{0}}\left( \mathbf{p^{2}} \rho^{-1} + \mathbf{q^{2}}\mathbf{J^{-1}}%
\right) dv,
\end{equation}
where 
\begin{equation}
\mathbf{p} = \frac{\partial \Upsilon _{C}}{\partial \mathbf{\dot{u}}}= \rho 
\frac{\partial \mathbf{u}}{\partial t} \text{ and } \mathbf{q} = \frac{%
\partial \Upsilon _{C}}{\partial \boldsymbol{\dot{\phi}}}= \mathbf{J} \frac{%
\partial \boldsymbol{\phi }}{\partial t},
\end{equation}
and 
\begin{equation}
\frac{\partial \mathbf{u}}{\partial t}=\frac{\partial \Upsilon _{K}}{%
\partial \mathbf{p}} = \mathbf{p}\rho^{-1} \text{ and } \frac{\partial %
\boldsymbol{\phi }}{\partial t}\mathbf{=}\frac{\partial \Upsilon _{K}}{%
\partial \mathbf{q}} = \mathbf{q J^{-1}},  \label{Const 1-1C}
\end{equation}

The work $T_{W}$ done by the inertia forces over displacement and
microrotation is given as in \cite{Steinberg2015} 
\begin{equation}
T_{W} = \int_{B_{0}} {\Upsilon_{W}} dv = \int_{B_{0}} \left( \frac{\partial 
\mathbf{p}}{\partial t} \mathbf{\cdot} \mathbf{u} + \frac{\partial \mathbf{q}%
}{\partial t} \mathbf{\cdot} \boldsymbol{\phi} \right) dv
\end{equation}

Keeping in mind that the variation of $\mathbf{p}$ $\mathbf{u}$, $\mathbf{q}$, $\boldsymbol{\phi}$, $\delta \mathbf{u}$ and $\delta \boldsymbol{\phi}$ is zero at $t_{0}$ and $t_{k}$ we can integrate by parts 
\begin{equation*}
\int_{t_{0}}^{t_{k}}T_{K}dt=\frac{1}{2}
\int_{B_{0}}\left( \mathbf{p\cdot u+q\cdot \boldsymbol{\phi}}\right)
dv|_{t_{0}}^{t_{k}}-\int_{t_{0}}^{t_{k}}T_{W}dt
\end{equation*}
\begin{equation*}
\delta \int_{t_{0}}^{t_{k}}T_{K}=-\delta \int_{t_{0}}^{t_{k}}T_{W}
\end{equation*}
or 
\begin{equation*}
\delta T_{C}=\text{ }\delta T_{K}=-\delta T_{W}
\end{equation*}
and therefore 
\begin{equation*}
\int_{t_{0}}^{t_{k}}\int_{B_{0}}\left( \left( \mathbf{p}\cdot \delta \left( 
\frac{\partial \mathbf{u}}{\partial t}\right) \mathbf{+q}\cdot \delta \left( 
{\frac{\partial \mathbf{\phi }}{\partial t}}\right) \right) \right)
dvdt = 
\end{equation*}

\begin{equation*}
-\int_{t_{0}}^{t_{k}}\int_{B_{0}}\left( \frac{\partial \mathbf{p}}{\partial t}\cdot \delta \mathbf{u+}\frac{\partial \mathbf{q}}{\partial t} \cdot \delta \boldsymbol{\phi }\right) dvdt  
\end{equation*}

\subsection{Variational Principle}

We modify the HPR principle \cite{Gurtin1972} for the case of Cosserat
elastodynamics in the following way. Now it states, that for any set $
\mathcal{A}$ of all admissible states $\mathfrak{s=}\left[ \mathbf{u},
\boldsymbol{\phi},\boldsymbol{\gamma},\boldsymbol{\chi},\boldsymbol{\sigma},
\boldsymbol{\mu}\right]$ that satisfy the {strain-displacement and
torsion-rotation relations (\ref{kinematic formulas}),} the zero variation 
\begin{equation*}
\delta \Theta (\mathfrak{s})=0
\end{equation*}%
of the functional 
\begin{eqnarray*}
\Theta (\mathfrak{s}) &=&
\int_{t_{0}}^{t_{k}}\left( U_{K}+T_{C}\right) dt + \\
&&\int_{t_{0}}^{t_{k}} \int_{B_{0}} {\left( \boldsymbol{\sigma}\cdot \boldsymbol{\gamma}+\boldsymbol{\mu}\cdot 
\boldsymbol{\chi}+\mathbf{p}\frac{\partial \mathbf{u}}{\partial t}\mathbf{q}{\frac{\partial \mathbf{\phi }}{\partial t}}\right)} dvdt +\\
&&\int_{t_{0}}^{t}\int_{\mathcal{G}_{1}}\left[ \boldsymbol{\sigma }_{%
\mathbf{n}}\cdot (\mathbf{u-u}_{\mathbf{0}})+\boldsymbol{\mu }_{\mathbf{n}%
}\left( \boldsymbol{\phi}-\boldsymbol{\phi}_{\mathbf{0}}\right) \right] dadt +\\
&&\int_{t_{0}}^{t}\int_{\mathcal{G}_{2}}\left[ \boldsymbol{\sigma }_{%
\mathbf{0}}\cdot \mathbf{u}+\boldsymbol{\mu}_{\mathbf{0}}\cdot %
\boldsymbol{\phi }\right] dadt  \notag
\end{eqnarray*}
at $\mathfrak{s\in }\mathcal{A}$ is equivalent of $\mathfrak{s}$ to be a solution of the system of equilibrium equations (\ref%
{equilibrium_equations 0}) - (\ref{equilibrium_equations}), constitutive relations (\ref{Hooke's_law 2}) - (\ref{Hooke's_law 2A}), which satisfies the mixed boundary conditions (\ref{Bound conditions}) - (\ref{Bound conditions a}). 

\hspace{0.1in}

\textbf{Proof of the Principle} \newline
Let us consider the variation of the functional $\Theta (\mathfrak{s})$:
\begin{eqnarray*}
\delta \Theta (\mathfrak{s}) &=& \int_{t_{0}}^{t_{k}}\left[ \delta U_{K}+\delta T_{C}\right] dt -\\
&&\int_{t_{0}}^{t_{k}}\int_{B_{0}} { \left(\delta \boldsymbol{\sigma}\cdot 
\boldsymbol{\gamma}+\boldsymbol{\sigma}\cdot \delta \boldsymbol{\gamma} \right) }dvdt +\\
&&\int_{t_{0}}^{t_{k}}\int_{B_{0}} { \left( \delta \boldsymbol{\mu}\cdot \boldsymbol{\chi}+\boldsymbol{\mu}\cdot \delta\boldsymbol{\chi} \right) }dvdt +\\
&&\int_{t_{0}}^{t_{k}}\int_{B_{0}} {\left[ \frac{\partial \mathbf{u}}{\partial t}\mathbf{\delta \mathbf{p}+p}\cdot \delta \left( \frac{\partial \mathbf{u}}{\partial t}\right) \right] }dvdt +\\
&&\int_{t_{0}}^{t_{k}}\int_{B_{0}} {\left[ \delta \mathbf{q}\cdot {\frac{\partial \mathbf{\phi}}{\partial t}}+\mathbf{q}\cdot \delta \left( {\frac{\partial \mathbf{\phi }}{\partial t}}\right) \right] }dvdt +\\
&&\int_{t_{0}}^{t}\int_{\mathcal{G}_{1}}\left( \delta \boldsymbol{
\sigma_{n}}\cdot (\mathbf{u-u}_{\mathbf{0}})+\boldsymbol{\sigma_{n}}\delta \boldsymbol{u}\right) dadt + \\
&&\int_{t_{0}}^{t}\int_{\mathcal{G}_{1}}\left(\delta \boldsymbol{\mu_{n}}\cdot \left( \boldsymbol{\phi}-\boldsymbol{\phi}_{\mathbf{0}}\right) +\boldsymbol{\mu_{n}}\delta \boldsymbol{\phi}\right) dadt + \\
&&\int_{t_{0}}^{t}\int_{\mathcal{G}_{2}}\left[ \boldsymbol{\sigma}_{\mathbf{0}}\cdot \delta \boldsymbol{u}+\boldsymbol{\mu}_{\mathbf{0}}\cdot \delta \boldsymbol{\phi}\right] dadt
\end{eqnarray*}
Taking into account (\ref{kinematic formulas}) we can perform the
integration by parts 
\begin{equation*}
\int_{B_{0}}\boldsymbol{\sigma}\cdot \delta \boldsymbol{\gamma}dv
=\int_{\partial B_{0}}\boldsymbol{\sigma_{n}}\cdot \delta \boldsymbol{u} da-\int_{B_{0}}{\left( \delta \boldsymbol{u}\cdot \mathbf{div}{\boldsymbol{\sigma}} - \boldsymbol{\varepsilon}\boldsymbol{\sigma}\cdot \delta \boldsymbol{\phi}\right)}dv\\
\end{equation*}%
and similar 
\begin{equation*}
\int_{B_{0}}\boldsymbol{\mu}\cdot \delta \boldsymbol{\chi}dv
=\int_{\partial B_{0}}\boldsymbol{\mu_{n}}\cdot \delta \boldsymbol{\phi}da-\int_{B_{0}}\delta \boldsymbol{\phi}\cdot \mathbf{div}{\boldsymbol{\mu}}dv
\end{equation*}%
Based on (\ref{Const1-1B}) - (\ref{Const 1-1C}) 
\begin{eqnarray}
\delta \Phi  &=&\frac{\partial \Phi }{\partial \boldsymbol{\sigma}}\cdot
\delta \boldsymbol{\sigma}+\frac{\partial \Phi }{\partial \boldsymbol{\mu}}%
\cdot \delta \boldsymbol{\mu},  \notag \\
\delta \Upsilon _{C} &=&\frac{\partial \Upsilon _{C}}{\partial \left( \frac{%
\partial \mathbf{u}}{\partial t}\right) }\cdot \delta \left( \frac{\partial 
\mathbf{u}}{\partial t}\right) +\frac{\partial \Upsilon _{C}}{\partial
\left( \frac{\partial \mathbf{\phi }}{\partial t}\right) }\cdot \delta
\left( \frac{\partial \mathbf{\phi }}{\partial t}\right) .  \notag
\end{eqnarray}
Then keeping in mind that $\delta T_{K}=-\delta T$ we can rewrite the expression for the variation of the functional $\delta \Theta (\mathfrak{s})$ in the following form 
\begin{eqnarray*}
\delta \Theta (\mathfrak{s}) &=&\int_{t_{0}}^{t}\int_{B_{0}}\left[ \left( 
\frac{\partial \Phi }{\partial \boldsymbol{\sigma}}-\boldsymbol{\gamma}%
\right) \cdot \delta \boldsymbol{\sigma}\right] dvdt+\\
&&\int_{t_{0}}^{t}\int_{B_{0}}\left[ \left( \frac{\partial \Phi }{\partial \boldsymbol{\mu}}-\boldsymbol{\chi}\right) \cdot \delta \boldsymbol{\mu}\right] dvdt + \notag \\
&&\int_{t_{0}}^{t}\int_{B_{0}}\left[ \left( \rho \frac{\partial \mathbf{u}}{%
\partial t}-\mathbf{p}\right) \cdot \delta \left( \frac{\partial \mathbf{u}}{%
\partial t}\right) \right] dvdt+\\
&&\int_{t_{0}}^{t}\int_{B_{0}}\left[ \left( 
\mathbf{J}\frac{\partial \mathbf{\varphi }}{\partial t}-\mathbf{q}\right)
\cdot \delta \left( \frac{\partial \mathbf{\varphi }}{\partial t}\right) %
\right] dvdt + \notag \\
&&\int_{t_{0}}^{t}\int_{B_{0}}\left[ \left( \mathbf{div}{\boldsymbol{\sigma}%
}-\frac{\partial \mathbf{p}}{\partial t}\right) \cdot \delta \boldsymbol{u}%
\right] dvdt+\\
&&\int_{t_{0}}^{t}\int_{B_{0}}\left[ \left( \mathbf{div}{%
\boldsymbol{\mu}}+\boldsymbol{\varepsilon}\cdot \boldsymbol{\sigma}-\frac{%
\partial \mathbf{q}}{\partial t}\right) \cdot \delta \boldsymbol{u}\right]
dvdt + \notag \\
&&\int_{t_{0}}^{t}\int_{\mathcal{G}_{1}}\left[ (\mathbf{u-u}_{\mathbf{0}%
})\cdot \delta \boldsymbol{\sigma_{n}}\right] dadt+\\
&&\int_{t_{0}}^{t}\int_{%
\mathcal{G}_{1}}\left[ (\boldsymbol{\phi-\phi}_{\mathbf{0}})\cdot \delta %
\boldsymbol{\mu_{n}}\right] dadt + \notag \\
&&\int_{t_{0}}^{t}\int_{\mathcal{G}_{2}}\left[ (\boldsymbol{\sigma}_{%
\mathbf{n}}-\boldsymbol{\sigma}_{\mathbf{0}})\cdot \delta \boldsymbol{u}%
\right] dadt+\\
&&\int_{t_{0}}^{t}\int_{\mathcal{G}_{2}}\left[ (\boldsymbol{\mu}_{%
\mathbf{n}}-\boldsymbol{\mu}_{\mathbf{0}})\cdot \delta \boldsymbol{\phi}%
\right] dadt  \notag
\end{eqnarray*}

\section{Dynamic Cosserat Plate Theory}

In this section we review our stress, couple stress and kinematic
assumptions of the Cosserat plate \cite{Steinberg2013}. We consider the thin plate $P$, where $h$ is the thickness of the plate and $x_{3}=0$ represent its middle plane. The sets $T$ and $B$ are the top and bottom surfaces contained in the planes $x_{3}=h/2$, $x_{3}=-h/2$ respectively and the curve  $\Gamma $ is the boundary of the middle plane of the plate.

The set of points $P=\left( \Gamma \times \lbrack -\frac{h}{2},\frac{h}{2}]\right) \cup T\cup B$ forms the entire surface of the plate and $\Gamma_{u}\times \lbrack -\frac{h}{2},\frac{h}{2}]$ is the lateral part of the boundary where displacements and microrotations are prescribed. The notation $\Gamma _{\sigma }=\Gamma \backslash \Gamma _{u}$ of the remainder we use to
describe the lateral part of the boundary edge $\Gamma _{\sigma }\times \lbrack -\frac{h}{2},\frac{h}{2}]$ where stress and couple stress are prescribed. We also use notation $P_{0}$ for the middle plane internal domain of the plate.

In our case we consider the vertical load and pure twisting momentum boundary conditions at the top and bottom of the plate, which can be written in the form:
\begin{eqnarray*}
\sigma_{33}(x_{1},x_{2},h/2,t)&=&\sigma ^{t}(x_{1},x_{2},t),\\\sigma_{33}(x_{1},x_{2},-h/2,t)&=&\sigma ^{b}(x_{1},x_{2},t),\\
\sigma_{3\beta }(x_{1},x_{2},\pm h/2,t) &=&0,\\
\mu_{33}(x_{1},x_{2},h/2,t) &=&\mu ^{t}(x_{1},x_{2},t),\\
\mu_{33}(x_{1},x_{2},-h/2,t)&=&\mu ^{b}(x_{1},x_{2},t),\\
\mu_{3\beta }(x_{1},x_{2},\pm h/2,t) &=&0,
\end{eqnarray*}
where $(x_{1},x_{2})\in P_{0}$.

We will also consider the rotatory inertia $\mathbf{J}$ in the form 
\begin{equation*}
\mathbf{J}=\left( 
\begin{matrix}
J_{11} & J_{12} & 0 \\ 
J_{12} & J_{22} & 0 \\ 
0 & 0 & J_{33}%
\end{matrix}%
\right) 
\end{equation*}

\subsection{Variational Principle for Dynamic Cosserat Plate}

Let $\mathcal{A}$ denote the set of all admissible states that satisfy the Cosserat plate strain-displacement relation (\ref{kinematic formulas}) and let $\Theta $ be a functional on $\mathcal{A}$ defined by
\begin{eqnarray*}
\Theta ({\large s,}\eta {\large )}&=&\int\limits_{P_{0}}\left(\mathcal{S\cdot E+}\mathbf{P}\cdot \frac{\partial \mathcal{U}}{\partial t}-%
\hat{P}\cdot \mathcal{W}+v\Omega _{3}^{0}\right)da+ \\
&&\int_{\Gamma _{\sigma }}\mathcal{S}_{o}\mathcal{\cdot }\left( \mathcal{U-U}_{o}\right)ds+ \\
&&\int_{\Gamma _{u}}\mathcal{S}_{n}\mathcal{\cdot U}ds +U_{K}^{S}+T_{C}^{S},  \end{eqnarray*}
for every ${\large s}=\left[ \mathcal{U},\mathcal{E},\mathcal{S}\right] \in \mathcal{A}.$ Here $\hat{P}=\left( \hat{p}_{1},\hat{p}_{2}\right)$, $\mathcal{W}=\left( W,W^{\ast }\right) $, $\widehat{p}_{1}=\eta p$, $\widehat{p}_{2}=\frac{2}{3}\left( 1-\eta \right) p$.

The plate stress and kinetic energy density are defined by the formulas 
\begin{eqnarray*}
U_{K}^{\mathcal{S}}=\int_{P_{0}}{\Phi \left( \mathcal{S} \right)}da, \\ 
T_{K}^{\mathcal{S}}=\int_{P_{0}}{\Upsilon_{C} \left( \frac{\partial \mathcal{U}}{\partial t} \right)} da,
\end{eqnarray*}
where $P_{0}$ is the internal domain of the middle plane of the plate and ${\Phi \left( \mathcal{S} \right)}$ and $\Upsilon_{C} \left( \frac{\partial \mathcal{U}}{\partial t} \right)$ are given as follows:
\begin{eqnarray*}
{\large \Phi (}\mathcal{S}{\large )} &=&-\frac{3\lambda \left(
M_{\alpha \alpha }\right) \left( M_{\beta \beta }\right) }{h^{3}\mu(3\lambda +2\mu )} +\frac{3\left( \alpha +\mu \right) M_{\alpha \beta }^{2}}{2h^{3}\alpha \mu} +\notag \\
&&\frac{3\left( \alpha +\mu \right) }{160h^{3}\alpha \mu }\left( 8\hat{Q}_{\alpha }\hat{Q}_{\alpha }+15Q_{\alpha }\hat{Q}_{\alpha }\right) +  \notag \\
&&\frac{3\left( \alpha +\mu \right) }{160h^{3}\alpha \mu }\left( 20\hat{Q}_{\alpha}Q_{\alpha }^{\ast }+8Q_{\alpha }^{\ast }Q_{\alpha }^{\ast }\right) + \notag \\
&& \frac{3\left( \gamma +\epsilon \right) S_{\alpha }^{\ast }S_{\alpha}^{\ast }}{2h^{3}\gamma \epsilon }+\frac{3\left( \alpha -\mu \right) M_{\alpha \beta }^{2}}{2h^{3}\alpha \mu 
}+\\
&&\frac{\alpha -\mu }{280h^{3}\alpha \mu }\left[ 21Q_{\alpha }\left( 5\hat{Q}_{\alpha }+4Q_{\alpha }^{\ast }\right) \right] - \notag \\
&&\frac{\gamma -\epsilon }{160h\gamma \epsilon }\left[ 24R_{\alpha \alpha
}^{2}+45R_{\alpha \alpha }^{\ast }+60R_{\alpha \beta }R_{\alpha \beta}^{\ast}\right] + \notag \\
&&\frac{\gamma +\epsilon }{160h^{3}\gamma
\epsilon }\left[ 8R_{\alpha \beta }^{2}+15R_{\alpha \beta }^{\ast }R_{\alpha\beta }^{\ast }+20R_{\alpha \beta }R_{\alpha \beta }^{\ast }\right] + \notag \\
&&\frac{3\beta }{80h\gamma (3\beta +2\gamma )}\left[8\left( R_{\alpha
\alpha }\right) \left( R_{\beta \beta }\right) +15\left( R_{\alpha \alpha}^{\ast }\right) \left( R_{\beta \beta }^{\ast }\right) \right]   + \notag \\
&&\frac{3\beta }{80h\gamma (3\beta +2\gamma )}\left[20\left( R_{\alpha
\alpha }\right) \left( R_{\alpha \alpha }^{\ast }\right) \right]   - \notag \\
&&\frac{\beta }{4\gamma (3\beta +2\gamma )}\left[ (2R_{\alpha \alpha}+3R_{\alpha \alpha }^{\ast })t-h\left( V^{2}+T^{2}\right) \right] + \notag
\\
&&\frac{\lambda }{560h\mu (3\lambda +2\mu )}\left[ \frac{5+3\eta }{(1+\eta )%
}pM_{\alpha \alpha }\right] + \notag \\
&&\frac{\left( \lambda +\mu \right) h}{840\mu (3\lambda +2\mu )}\left(\frac{140+168\eta +51\eta ^{2}}{4(1+\eta )^{2}}\right) p^{2}+ \notag \\
&&\frac{\left(\lambda +\mu \right) h}{2\mu (3\lambda +2\mu )}\sigma _{0}^{2}+\frac{\epsilon h}{12h\gamma \epsilon }\left[ (3T^{2}+V^{2})\right] 
\end{eqnarray*}
and 
\begin{eqnarray*}
{\LARGE \Upsilon }_{C}\left( \frac{\partial \mathcal{U}}{\partial t}\right) 
&=&\frac{h\rho}{2}\left( \frac{\partial W}{\partial t}\right) ^{2}+\frac{4h\rho }{15}\left( \frac{\partial W^{\ast }}{\partial t}\right) ^{2}+\\
&&\frac{2h\rho}{3}\left(\frac{\partial W}{\partial t}\frac{\partial W^{\ast }}{\partial t}\right) +\frac{h^{3}\rho }{24}\left( \frac{\partial \Psi _{\alpha
}}{\partial t}\right) ^{2}+ \\
&&\frac{4hJ_{\alpha \beta }}{15}\left( \frac{\partial \Omega _{\alpha }^{0}}{%
\partial t}\frac{\partial \Omega _{\beta }^{0}}{\partial t}\right) +\frac{hJ_{\alpha \beta }}{2}\left( \frac{\partial \hat{\Omega}_{\alpha }^{0}}{\partial t}\frac{\partial \hat{\Omega}_{\beta }^{0}}{\partial t}\right)+
\\ &&\frac{2hJ_{\alpha \beta }}{3}\left( \frac{\partial \Omega _{\alpha }^{0}}{%
\partial t}\frac{\partial \hat{\Omega}_{\beta }^{0}}{\partial t}\right) +%
\frac{hJ_{33}}{6}\left( \frac{\partial \Omega _{3}}{\partial t}\right) ^{2}.
\end{eqnarray*}

$\mathcal{S}$, $\mathcal{U}$ and $\mathcal{E}$ are the Cosserat plate stress, displacement and strain sets 
\begin{eqnarray}
\mathcal{S} &=&\left[ M_{\alpha \beta },Q_{\alpha },Q_{\alpha }^{\ast },\hat{%
Q}_{\alpha },R_{\alpha \beta },R_{\alpha \beta }^{\ast },S_{\beta }^{\ast }%
\right] , \\
\mathcal{S}_{n} &=&\left[ \check{M}_{\alpha },\check{Q}^{\ast },\check{Q}%
_{\alpha }^{\symbol{94}},\check{R}_{\alpha },\check{R}_{\alpha }^{\ast },%
\check{S}^{\ast }\right] , \\
\mathcal{S}_{o} &=&\left[ \Pi _{o\alpha },\Pi _{o3},\Pi _{o3}^{\ast
},M_{o\alpha },M_{o\alpha }^{\ast },M_{o3}^{\ast }\right] , \\
\mathcal{U} &=&\left[ \Psi _{\alpha },W,\Omega _{3},\Omega _{\alpha
}^{0},W^{\ast },\Omega _{\alpha }^{0}\right] , \\
\mathcal{E} &=&\left[ e_{\alpha \beta },\omega _{\beta },\omega _{\alpha
}^{\ast },\hat{\omega}_{\alpha },\tau _{3\alpha },\tau _{\alpha \beta },\tau
_{\alpha \beta }^{\ast }\right] ,
\end{eqnarray}
where 
\begin{eqnarray*}
M_{\alpha \beta }n_{\beta}=\Pi _{o\alpha }, R_{\alpha \beta }n_{\beta}=M_{o\alpha }, \\
Q_{\alpha }^{\ast }n_{\alpha} = \Pi_{o3}, S_{\alpha }^{\ast }n_{\alpha}=M_{o3}^{\ast}, \\
\hat{Q}_{\alpha }n_{\alpha}=\Pi _{o3}^{\ast},R_{\alpha \beta }^{\ast}n_{\beta }=M_{o\alpha }^{\ast },\\
\check{M}_{\alpha }=M_{\alpha \beta }n_{\beta}, \check{Q}^{\ast}=Q_{\beta }^{\ast }n_{\beta},\\
\check{R}_{\alpha }=R_{\alpha \beta}n_{\beta}, \check{S}^{\ast} =S_{\beta }^{\ast }n_{\beta},\\
\check{Q}^{\symbol{94}}=\check{Q}_{\beta }^{\symbol{94}}n_{\beta }, \check{R}_{\alpha}^{\ast }=\check{R}_{\alpha \beta}^{\ast}n_{\beta}.
\end{eqnarray*}
In the above $n_{\beta}$ is the outward unit normal vector to $\Gamma_{u}$.

The plate characteristics, being the functions of $x_{1}$, $x_{2}$ and $t$, provide the approximation of the components of the three-dimensional tensors $\boldsymbol{\sigma}$ and $\boldsymbol{\mu}$
\begin{eqnarray}
\sigma_{\alpha \beta } &=& \frac{6}{h^{2}}\zeta M_{\alpha \beta
},  \label{stress assumption 1} \\
\sigma _{3\beta } &=& \frac{3}{2h}\left( 1-\zeta ^{2}\right) Q_{\beta}, \label{stress assumption 1a} \\
\sigma_{\beta 3} &=& \frac{3}{2h}\left( 1-\zeta ^{2}\right) Q_{\beta }^{\ast}+\frac{3}{2h}\hat{Q}_{\beta }, \label{stress assumption 1b} \\
\sigma_{33} &=& -\frac{3}{4}\left( \frac{1}{3}\zeta ^{3}-\zeta \right) p_{1}+\zeta p_{2}+\sigma _{0}, \label{stress assumption 1ca} \\
\mu _{\alpha \beta } &=& \frac{3}{2h}\left( 1-\zeta ^{2}\right) R_{\alpha \beta }+\frac{3}{2h}R_{\alpha \beta }^{\ast }, \label{stress assumption 1d} \\
\mu _{\beta 3} &=& \frac{6}{h^{2}}\zeta S_{\beta }^{\ast}, \label{couple stress 1} \\
\mu _{3\beta } &=& 0,  \label{couple stress 1a} \\
\mu _{33} &=& \zeta V+T, \label{stress assump 3(a)}
\end{eqnarray}
where 
\begin{eqnarray}
p &=&\sigma ^{t}(x_{1},x_{2},t)-\sigma ^{b}(x_{1},x_{2},t), \\
\sigma _{0} &=&\frac{1}{2}\left( \sigma
^{t}(x_{1},x_{2},t)+\sigma ^{b}(x_{1},x_{2},t)\right) , \\
V &=&\frac{1}{2}\left( \mu ^{t}(x_{1},x_{2},t)-\mu
^{b}(x_{1},x_{2},t)\right) , \\
T &=&\frac{1}{2}\left( \mu ^{t}(x_{1},x_{2},t)+\mu
^{b}(x_{1},x_{2},t)\right) , \\
p_{1} &=&\eta p(x_{1},x_{2},t), \\
p_{2} &=&\frac{\left( 1-\eta \right) }{2}p(x_{1},x_{2},t),
\end{eqnarray}
three-dimensional displacements $\mathbf{u}$ and microrotations $\boldsymbol{\phi}$  
\begin{eqnarray}
u_{\alpha} &=& \frac{h}{2}\zeta \Psi_{\alpha}, \\
u_{3} &=& W+\left( 1-\zeta^{2} \right) W^{\ast}, \\
\phi_{\alpha} &=& \Omega_{\alpha}^{0} \left( 1-\zeta^{2} \right) +\hat{\Omega}_{\alpha}, \\
\phi_{3} &=&\zeta \Omega_{3},
\end{eqnarray}
and the three-dimensional strain and torsion tensors $\mathbf{\gamma}$ and $\mathbf{\chi}$ 
\begin{eqnarray}
\gamma _{\alpha \beta } &=&\frac{6}{h^{2}}\zeta e_{\alpha \beta
},  \label{strain assumption 1} \\
\gamma _{3\beta } &=&\frac{3}{2h}\left( 1-\zeta ^{2}\right) \omega _{\beta},  \label{strain assumption 1a} \\
\gamma _{\beta 3} &=&\frac{3}{2h}\left( 1-\zeta ^{2}\right) \omega _{\beta
}^{\ast }+\frac{3}{2h}\hat{\omega}_{\beta },
\label{strain assumption 1ca} \\
\chi _{\alpha \beta } &=&\frac{3}{2h}\left( 1-\zeta ^{2}\right) \tau_{\alpha \beta }+\frac{3}{2h}\tau _{\alpha \beta }^{\ast},  \label{couple strain 1} \\
\chi _{3\beta } &=&\frac{6}{h^{2}}\zeta \tau _{\beta }^{\ast},  \label{couple strain 1a}
\end{eqnarray}
where $\zeta =\frac{2x_{3}}{h}$.

Then zero variation of the functional  
\begin{equation*}
\delta \Theta ({\large s,}\eta )=0
\end{equation*}%
is equivalent to the plate bending system of equations (A) and constitutive formulas (B) mixed problems.

(A). The bending equilibrium system of equations: 
\begin{eqnarray}
M_{\alpha \beta ,\alpha }-Q_{\beta } &=&I_{1}\frac{\partial ^{2}\Psi _{\beta
}}{\partial t^{2}},  \label{equilibrium_equations 1_E} \\
Q_{\alpha ,\alpha }^{\ast }+\hat{p}_{1} &=&I_{2}\frac{\partial ^{2}W^{\ast }%
}{\partial t^{2}},  \label{equilibrium_equations 1_B} \\
R_{\alpha \beta ,\alpha }+\varepsilon _{3\beta \gamma }\left( Q_{\gamma
}^{\ast }-Q_{\gamma }\right)  &=&I_{\alpha \beta }\frac{\partial ^{2}\Omega
_{\alpha }^{0}}{\partial t^{2}}, \\
\varepsilon _{3\beta \gamma }M_{\beta \gamma }+S_{\alpha ,\alpha }^{\ast }
&=&I_{3}\frac{\partial ^{2}\Omega _{3}}{\partial t^{2}},
\label{equilibrium_equations 1_C} \\
\hat{Q}_{\alpha ,\alpha }+\hat{p}_{2} &=&I_{2}\frac{\partial ^{2}W}{\partial
t^{2}},  \label{equilibrium_equations 1_DD} \\
R_{\alpha \beta ,\alpha }^{\ast }+\varepsilon _{3\beta \gamma }\hat{Q}%
_{\gamma } &=&I_{\alpha \beta }^{0}\frac{\partial ^{2}\hat{\Omega}_{\alpha
}^{0}}{\partial t^{2}},  \label{equilibrium_equations 1_D}
\end{eqnarray}%
where $I_{1}=\frac{h^{3}}{12}\rho $, $I_{2}=\frac{2h}{3}\rho $, $I_{\alpha
\beta }=\frac{5h}{6}J_{\alpha \beta }$, $I_{3}=\frac{h^{2}}{6}J_{33}$, $%
I_{\alpha \beta }^{0}=\frac{2h}{3}J_{\alpha \beta }$, $\widehat{p}_{1}=\eta
_{opt}p$ and $\widehat{p}_{2}=\frac{2}{3}\left( 1-\eta _{opt}\right) p$,
with the resultant traction boundary conditions : 
\begin{eqnarray}
M_{\alpha \beta }n_{\beta } = \Pi _{o\alpha },\ R_{\alpha \beta }n_{\beta
}=M_{o\alpha },  \label{bc_0} \\
Q_{\alpha }^{\ast }n_{\alpha } = \Pi _{o3},\ S_{\alpha }^{\ast }n_{\alpha
}=\Upsilon _{o3},  \label{bc_1}
\end{eqnarray}%
at the part $\Gamma _{\sigma }$ and the resultant displacement boundary
conditions

\begin{equation}
\Psi _{\alpha }=\Psi _{o\alpha },\text{ }W=W_{o},\text{ }\Omega _{\alpha
}^{0}=\Omega _{o\alpha }^{0},\text{ }\Omega _{3}=\Omega _{o3},  \label{bu_1}
\end{equation}%
at the part $\Gamma _{u}.$

(B). Constitutive formulas in the reverse form : $\footnote{%
In the following formulas a subindex $\beta =1$ if $\alpha =2$ and $\beta =2$
if $\alpha =1$}$

\begin{eqnarray*}
M_{\alpha \alpha }=\frac{\mu \left( \lambda +\mu \right) h^{3}}{3(\lambda
+2\mu )}\Psi _{\alpha ,\alpha }+\frac{\lambda \mu h^{3}}{6(\lambda +2\mu )}%
\Psi _{\beta ,\beta }+\frac{\left( 3p_{1}+5p_{2}\right) \lambda h^{2}}{%
30(\lambda +2\mu )}, \\
M_{\beta \alpha }=\frac{\left( \mu -\alpha \right) h^{3}}{12}\Psi _{\alpha
,\beta }+\frac{\left( \mu +\alpha \right) h^{3}}{12}\Psi _{\beta ,\alpha
}+(-1)^{\alpha ^{\prime }}\frac{\alpha h^{3}}{6}\Omega _{3},  \end{eqnarray*}
\begin{eqnarray*}
R_{\beta \alpha } = \frac{5\left( \gamma -\epsilon \right) h}{6}\Omega
_{\beta ,\alpha }^{0}+\frac{5\left( \gamma +\epsilon \right) h}{6}\Omega
_{\alpha ,\beta }^{0}, \\
R_{\alpha \alpha } = \frac{10h\gamma (\beta +\gamma )}{3(\beta +2\gamma )}%
\Omega _{\alpha ,\alpha }^{0}+\frac{5h\beta \gamma }{3(\beta +2\gamma )}%
\Omega _{\beta ,\beta }^{0}, \\
R_{\beta \alpha }^{\ast } = \frac{2\left( \gamma -\epsilon \right) h}{3}\hat{\Omega}%
_{\beta ,\alpha }+\frac{2\left( \gamma +\epsilon \right) h}{3}\hat{\Omega}%
_{\alpha ,\beta }, \\
R_{\alpha \alpha }^{\ast } = \frac{8\gamma \left( \gamma +\beta \right) h}{3(\beta
+2\gamma )}\hat{\Omega}_{\alpha ,\alpha }+\frac{4\gamma \beta h}{3(\beta
+2\gamma )}\hat{\Omega}_{\beta ,\beta },
\end{eqnarray*}

\begin{eqnarray*}
Q_{\alpha } &=&\frac{5\left( \mu +\alpha \right) h}{6}\Psi _{\alpha }+\frac{%
5\left( \mu -\alpha \right) h}{6}W_{,\alpha }+\frac{2\left( \mu -\alpha
\right) h}{3}W_{,\alpha }^{\ast }+\\
&&(-1)^{\beta }\frac{5h\alpha }{3}\Omega
_{\beta }^{0}+(-1)^{\beta }\frac{5h\alpha }{3}\hat{\Omega}_{\beta }, \\
Q_{\alpha }^{\ast } &=&\frac{5\left( \mu -\alpha \right) h}{6}\Psi _{\alpha
}+\frac{5\left( \mu -\alpha \right) ^{2}h}{6\left( \mu +\alpha \right) }%
W_{,\alpha }+\frac{2\left( \mu +\alpha \right) h}{3}W_{,\alpha }^{\ast
}+\\
&&(-1)^{\alpha }\frac{5h\alpha }{3}\left( \Omega _{\beta }^{0}+\frac{\left(
\mu -\alpha \right) }{\left( \mu +\alpha \right) }\hat{\Omega}_{\beta
}\right) ,\\
\hat{Q}_{\alpha } &=&\frac{8\alpha \mu h}{3\left( \mu +\alpha \right) }%
W_{,\alpha }+(-1)^{\alpha }\frac{8\alpha \mu h}{3\left( \mu +\alpha \right) }
\hat{\Omega}_{\beta },
\end{eqnarray*}
\begin{equation}
S_{\alpha }^{\ast }=\frac{5\gamma \epsilon h^{3}}{3\left( \gamma +\epsilon
\right) }\Omega _{3,\alpha },  \label{Cons 5}
\end{equation}
and the optimal value $\eta _{opt}$ of the splitting parameter is given as
in \cite{Kvasov2013} 
\begin{equation}
\eta _{opt}=\frac{2\mathcal{W}^{\left( 00\right) }-\mathcal{W}^{\left(
10\right) }-\mathcal{W}^{\left( 01\right) }}{2\left( \mathcal{W}^{\left(
11\right) }+\mathcal{W}^{\left( 00\right) }-\mathcal{W}^{\left( 10\right) }-%
\mathcal{W}^{\left( 01\right) }\right) }.  \label{eta_minimum}
\end{equation}%
where $\mathcal{W}^{(ij)}=\left.\mathcal{S}\right\vert_{\eta=i} \cdot \left.\mathcal{E}\right\vert_{\eta=j}.$

We also assume that the initial condition can be presented in the  form 
\begin{eqnarray*}
\mathcal{U}\left( x_{1},x_{2},0\right) =\mathcal{U}^{0}\left(
x_{1},x_{2}\right) , \\
\frac{\partial \mathcal{U}}{\partial t}\left(
x_{1},x_{2},0\right) =\mathcal{V}^{0}\left( x_{1},x_{2}\right) 
\end{eqnarray*}

\subsection{Cosserat Plate Dynamic Field Equations}

The Cosserat plate field equations are obtained by substituting the relations into the system of equations (\ref{equilibrium_equations 1_E}) -- (\ref{equilibrium_equations 1_D})
similar to \cite{Kvasov2013}: 
\begin{equation}
L\mathcal{U}=K\frac{\partial ^{2}\mathcal{U}}{\partial t^{2}}+\mathcal{F}%
\left( \eta \right)  \label{bending_system}
\end{equation}%
where 
\begin{equation*}
L= \left[
\begin{array}{ccccccccc}
L_{11} & L_{12} & L_{13} & L_{14} & 0 & L_{16} & kL_{13} & 0 & L_{16} \\ 
L_{12} & L_{22} & L_{23} & L_{24} & L_{16} & 0 & kL_{23} & L_{16} & 0 \\ 
-L_{13} & -L_{23} & L_{33} & 0 & L_{35} & L_{36} & L_{77} & L_{38} & L_{39}
\\ 
L_{41} & L_{42} & 0 & L_{44} & 0 & 0 & 0 & 0 & 0 \\ 
0 & -L_{16} & -L_{38} & 0 & L_{55} & L_{56} & -kL_{35} & L_{58} & 0 \\ 
L_{16} & 0 & -L_{39} & 0 & L_{56} & L_{66} & -kL_{36} & 0 & L_{58} \\ 
-L_{13} & -L_{14} & L_{73} & 0 & L_{35} & L_{36} & L_{77} & L_{78} & L_{79}
\\ 
0 & -L_{16} & -L_{78} & 0 & L_{85} & L_{56} & -kL_{35} & L_{88} & kL_{56} \\ 
L_{16} & 0 & -L_{79} & 0 & L_{56} & L_{55} & -kL_{36} & kL_{56} & L_{99}%
\end{array} \right]
\end{equation*}

\begin{equation*}
K=
\left[
\begin{array}{ccccccccc}
\frac{h^{3}\rho}{12} & 0 & 0 & 0 & 0 & 0 & 0 & 0 & 0 \\ 
0 & \frac{h^{3}\rho}{12} & 0 & 0 & 0 & 0 & 0 & 0 & 0 \\ 
0 & 0 & \frac{2h\rho}{3} & 0 & 0 & 0 & 0 & 0 & 0 \\ 
0 & 0 & 0 & \frac{h^{2}J_{33}}{6} & 0 & 0 & 0 & 0 & 0 \\ 
0 & 0 & 0 & 0 & \frac{5hJ_{11}}{6} & \frac{5hJ_{12}}{6} & 0 & 0 & 0 \\ 
0 & 0 & 0 & 0 & \frac{5hJ_{12}}{6} & \frac{5hJ_{22}}{6} & 0 & 0 & 0 \\ 
0 & 0 & 0 & 0 & 0 & 0 & \frac{2h\rho}{3} & 0 & 0 \\ 
0 & 0 & 0 & 0 & 0 & 0 & 0 & \frac{2hJ_{11}}{3} & \frac{2hJ_{12}}{3} \\ 
0 & 0 & 0 & 0 & 0 & 0 & 0 & \frac{2hJ_{12}}{3} & \frac{2hJ_{22}}{3}%
\end{array} \right]
\end{equation*}

\begin{equation*}
\mathcal{F}\left( \eta \right) =
\left[
\begin{array}{ccccccccc}
-\frac{3h^{2}\lambda \left( 3p_{1,1}+5p_{2,1}\right) }{30\left( \lambda+2\mu \right) } \\ 
-\frac{3h^{2}\lambda \left(3p_{1,2}+5p_{2,2}\right) }{30\left( \lambda +2\mu \right) } \\ 
-p_{1} \\
0 \\
0 \\
0 \\
\frac{h^{2}(3p_{1}+4p_{2})}{24} \\
0 \\
0 
\end{array} \right]
\end{equation*}

Here $p_{1}=\eta p, p_{2}=\frac{\left( 1-\eta \right) }{2}p$ and $\mathcal{U}$ is given as before
\begin{equation*}
\mathcal{U}=
\begingroup 
\setlength\arraycolsep{2pt}
\begin{pmatrix}
\Psi _{1}, & \Psi _{2}, & W, & \Omega _{3}, & \Omega _{1}^{0}, & \Omega
_{2}^{0}, & W^{\ast }, & \Omega _{1}^{0}, & \Omega _{2}^{0}%
\end{pmatrix}
\endgroup ^{T}
\end{equation*}

The operators $L_{ij}$ are given as follows
\begin{align*}
L_{11}&=c_{1}\frac{\partial ^{2}}{\partial x_{1}^{2}}+c_{2}\frac{\partial^{2}}{\partial x_{2}^{2}}-c_{3}, & L_{12}&=(c_{1}-c_{2})\frac{\partial ^{2}}{\partial x_{1}x_{2}}, \\ 
L_{13}&=c_{11}\frac{\partial }{\partial x_{1}}, &
L_{14}&=c_{12}\frac{\partial }{\partial x_{2}}, \\ 
L_{16}&=c_{13}, & 
L_{17}&=k_{1}c_{11}\frac{\partial }{\partial x_{1}}, \\
L_{22}&=c_{2}\frac{\partial ^{2}}{\partial x_{1}^{2}}+c_{1}\frac{\partial^{2}}{\partial x_{2}^{2}}-c_{3}, & L_{23}&=c_{11}\frac{\partial}{\partial x_{2}}, \\
L_{24}&=-c_{12}\frac{\partial}{\partial x_{1}}, &
L_{33}&=c_{3}(\frac{\partial^{2}}{\partial x_{1}^{2}}+\frac{\partial^{2}}{\partial x_{2}^{2}}), \\
L_{35}&=-c_{13}\frac{\partial}{\partial x_{2}}, & 
L_{36}&=c_{13}\frac{\partial}{\partial x_{1}}, \\
L_{38}&=-c_{10}\frac{\partial }{\partial x_{2}}, & L_{39}&=c_{10}\frac{\partial }{\partial x_{1}}, \\
L_{41}&=-c_{12}\frac{\partial}{\partial x_{2}}, &
L_{42}&=c_{12}\frac{\partial }{\partial x_{1}}, \\
L_{44}&=c_{6}\left( \frac{\partial ^{2}}{\partial x_{1}^{2}}+\frac{\partial^{2}}{\partial x_{2}^{2}}%
\right) -2c_{12}, & L_{55}&=c_{7}\frac{\partial ^{2}}{\partial x_{1}^{2}}%
+c_{8}\frac{\partial ^{2}}{\partial x_{2}^{2}}-2c_{13}, \\
L_{56}&=(c_{7}-c_{8})\frac{\partial ^{2}}{\partial x_{1}x_{2}}, & 
L_{58}&=-c_{9}, \\
L_{66}&=c_{8}\frac{\partial ^{2}}{\partial x_{1}^{2}}+c_{7} \frac{\partial ^{2}}{\partial x_{2}^{2}}-2c_{13}, &
L_{73}&=c_{5}(\frac{\partial^{2}}{\partial x_{1}^{2}}+\frac{\partial ^{2}}{\partial x_{2}^{2}}),\\ 
L_{77}&=c_{4}(\frac{\partial ^{2}}{\partial x_{1}^{2}%
}+\frac{\partial ^{2}}{\partial x_{2}^{2}}), & L_{78}&=-c_{14}\frac{\partial%
}{\partial x_{2}}, \\
L_{79}&=c_{14}\frac{\partial}{\partial x_{1}}, & L_{85}&=c_{7}\frac{\partial^{2}}{\partial x_{1}^{2}}+c_{8}\frac{\partial ^{2}}{\partial x_{2}^{2}}-2c_{13}, \\
L_{88}&=c_{7}\frac{\partial ^{2}}{\partial x_{1}^{2}}+c_{8}\frac{\partial ^{2}}{\partial x_{2}^{2}}-c_{15}, &
L_{99}&=c_{8}\frac{\partial ^{2}}{\partial x_{1}^{2}}+c_{7}\frac{\partial^{2}}{\partial x_{2}^{2}}-c_{15}.
\end{align*}

The coefficients $c_{i}$ are given as 
\begin{align*}
c_{1}&=\frac{h^{3}\mu (\lambda +\mu )}{3(\lambda +2\mu )}, & c_{2}&=\frac{%
h^{3}(\alpha +\mu)}{12}, \\ c_{3}&=\frac{5h(\alpha +\mu )}{6}, & c_{4}&=\frac{%
5h(\alpha-\mu)^{2}}{6(\alpha+\mu )}, \\
c_{5}&=\frac{h(5\alpha ^{2}+6\alpha \mu+5\mu ^{2})}{6(\alpha +\mu )}, & 
c_{6}&=\frac{h^{3}\gamma \epsilon }{3(\gamma +\epsilon )}, \\
c_{7}&=\frac{10h\gamma \left( \beta +\gamma \right) }{3\left(\beta +2\gamma \right)}, & 
c_{8}&=\frac{5h\left( \gamma +\epsilon \right)}{6}, \\
c_{9}&=\frac{10h\alpha^{2}}{3(\alpha +\mu )}, & c_{10}&=\frac{5h\alpha
(\alpha-\mu )}{3(\alpha +\mu )},\\
c_{11}&=\frac{5h(\alpha -\mu )}{6}, & 
c_{12}&=\frac{h^{3}\alpha }{6}, \\
c_{13}&=\frac{5h\alpha }{3}, & c_{14}&=\frac{h\alpha (5\alpha+3\mu )}{3(\alpha +\mu )}, \\
c_{15}&=\frac{2h\alpha (5\alpha+4\mu )}{3(\alpha+\mu )}.& 
\end{align*}

\section{Numerical Simulation}

\subsection{Cosserat Plate Vibration}

In our computations we consider the plates made of polyurethane foam -- a material reported in the literature to behave Cosserat like -- and the values of the technical elastic parameters presented in \cite{Lakes1995}: 
\begin{eqnarray*}
E &=& 299.5 \text{ MPa}, \\
\nu &=& 0.44, \\
l_{t} &=& 0.62 \text{ mm}, \\
l_{b} &=& 0.327 \text{ mm}, \\
N^{2} &=& 0.04.
\end{eqnarray*}

Taking into account that the ratio $\beta /\gamma $ is equal to 1 for bending \cite{Lakes1995}, these values of the technical constants correspond to the following values of Lam\'{e} and Cosserat parameters: 
\begin{eqnarray*}
\lambda &=& 762.616 \text{ MPa}, \\
\mu &=& 103.993 \text{ MPa}, \\
\alpha &=& 4.333 \text{ MPa}, \\
\beta &=& 39.975 \text{ MPa}, \\
\gamma &=& 39.975 \text{ MPa}, \\
\epsilon &=& 4.505 \text{ MPa}.
\end{eqnarray*}

We consider a low density rigid foam usually characterized by the densities of 24-50 kg/m$^3$ \cite{Singh2002}. In all further numerical computations we used the density value $\rho=34$ kg/m$^3$ and different values the rotatory inertia $\mathbf{J}$.

\noindent
\begin{center}
\begin{table}[H]
\caption{Eigenfrequencies $\protect\omega _{i}^{11}$ (Hz) for
different shapes of micro-elements}
\label{tab:different}%
\setkeys{Gin}{keepaspectratio} {
\begin{tabularx}{\textwidth}{@{}p{.20\textwidth}XXXXXXXX@{}}
\hline
Shape & $J_{x}$ & $J_{y}$ & $J_{z}$ & $\omega_{1}$, $\omega_{2}$ & $\omega_{3}$, $\omega_{7}$ & $\omega_{4}$ & $\omega_{5}$, $\omega_{8}$ & $\omega_{6}$, $\omega_{9}$ \\ 
\hline
Ball & 0.001 & 0.001 & 0.001 & 17.88 & 0.31  & 501.13 & 205.62 & 338.95 \\
Vertical Ellipsoid & 0.001 & 0.001 & 0.0001 & 17.88 & 0.31 & 501.13 & 650.22 & 338.95 \\
Horizontal Ellipsoid & 0.0001 & 0.001 & 0.001 & 17.88 & 0.31  & 1363.01 & 205.62 & 394.08 \\
\hline
\end{tabularx}
}
\end{table}
\end{center}

We consider a plate $a\times a$ of thickness $h$ with the boundary $G =
G_{1}\cup G_{2}$ 
\begin{eqnarray*}
G_{1} &=&\left\{ \left( x_{1},x_{2}\right) :x_{1}\in \left\{ 0,a\right\}
,x_{2}\in \left[ 0,a\right] \right\} \\
G_{2} &=&\left\{ \left( x_{1},x_{2}\right) :x_{2}\in \left\{ 0,a\right\}
,x_{1}\in \left[ 0,a\right] \right\}
\end{eqnarray*}
and the following hard simply supported boundary conditions \cite%
{Steinberg2013}: 
\begin{eqnarray}
G_{1} &:&W=0,\text{ }W^{\ast }=0,\text{ }\Psi _{2}=0,\text{ }\Omega
_{1}^{0}=0,\text{ }\hat{\Omega}_{1}^{0}=0;\text{ } \\
G_{1} &:&\Omega _{3}=0,\text{ }\frac{\partial \Psi _{1}}{\partial n}=0,\text{
}\frac{\partial \Omega _{2}^{0}}{\partial n}=0,\text{ }\frac{\partial \hat{%
\Omega}_{2}^{0}}{\partial n}=0; \\
G_{2} &:&W=0,\text{ }W^{\ast }=0,\text{ }\Psi _{1}=0,\text{ }\Omega
_{2}^{0}=0,\text{ }\hat{\Omega}_{2}^{0}=0;\text{ } \\
G_{2} &:&\Omega _{3}=0,\text{ }\frac{\partial \Psi _{2}}{\partial n}=0,\text{
}\frac{\partial \Omega _{1}^{0}}{\partial n}=0,\text{ }\frac{\partial \hat{%
\Omega}_{1}^{0}}{\partial n}=0.
\end{eqnarray}

Similar to \cite{Steinberg2015} we apply the method of separation of
variables for the eigenvalue problem (\ref{bending_system}) to solve for the kinematic variables $\Psi_{\alpha}$, $W$, $\Omega_{3}$, $\Omega_{\alpha}^{0}$%
, $W^{\ast}$ and $\Omega_{\alpha}^{0}$. The kinematic variables can be
further expressed in the following form 
\begin{eqnarray}
\Psi_{1}^{nm} &=& \Psi_{1(A)}^{nm}  + \Psi_{1(B)}^{nm}, \label{psi1} \\
\Psi_{2}^{nm}&=& \Psi_{2(A)}^{nm}  + \Psi_{2(B)}^{nm},   \label{psi2} \\
W^{nm} &=& W_{(A)}^{nm} + W_{(B)}^{nm},  \label{w} \\
\Omega_{3}^{nm}&=& \Omega_{3(A)}^{nm}+\Omega_{3(B)}^{nm},  \label{omega3} \\
\Omega_{1}^{0,nm} &=& \Omega_{1(A)}^{0,nm}+\Omega_{1(B)}^{0,nm},  \label{omega1} \\
\Omega_{2}^{0,nm}&=& \Omega_{2(A)}^{0,nm}+\Omega_{2(B)}^{0,nm},  \label{omega2} \\
W^{\ast ,nm} &=& W_{(A)}^{\ast ,nm}+W_{(B)}^{\ast ,nm},  \label{ws} \\
\hat{\Omega}_{1}^{0,nm}&=& \hat{\Omega}_{1(A)}^{0,nm} + \hat{\Omega}_{1(B)}^{0,nm},  \label{omega01} \\
\hat{\Omega}_{2}^{0,nm} &=& \hat{\Omega}_{2(A)}^{0,nm} + \hat{\Omega}_{2(B)}^{0,nm},  \label{omega02}
\end{eqnarray}

\noindent 
\begin{figure}[H]
\begin{center}
\includegraphics[width=3.7in]{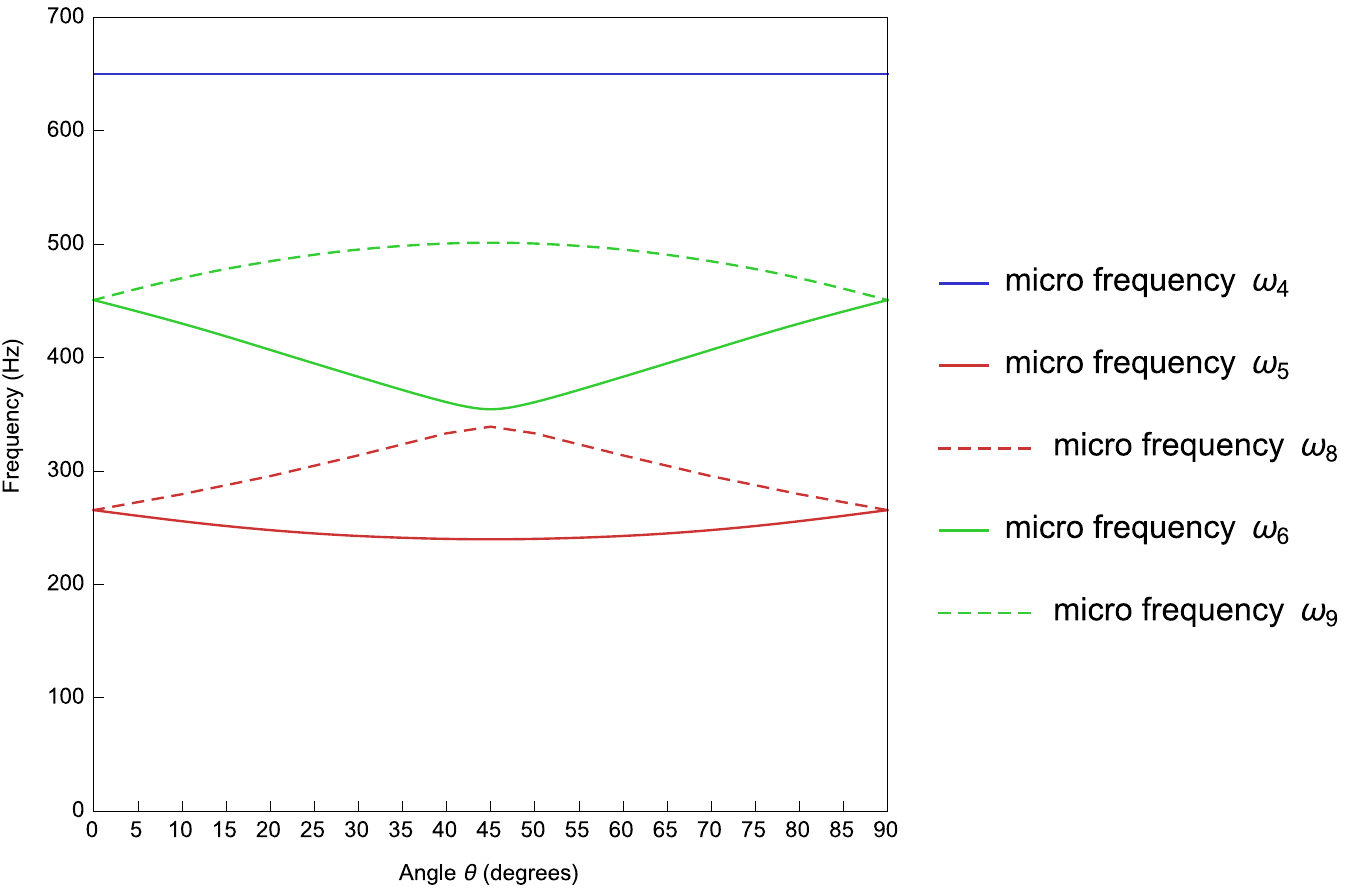}
\end{center}
\caption{\textit{Micro frequencies $\omega_{4}$,  $\omega_{5}$,  $\omega_{8}$,  $\omega_{6}$ and  $\omega_{9}$.}}
\label{fig:microfreq}
\end{figure}

where
\begin{eqnarray*}
\Psi_{1(A)}^{nm} &=& A_{1}\cos \left( \frac{n\pi x_{1}}{a}\right) \sin \left(\frac{m\pi x_{2}}{a}\right)\sin \left( \omega t\right),  \label{psi1} \\
\Psi_{2(A)}^{nm} &=& A_{2}\sin \left( \frac{n\pi x_{1}}{a}\right) \cos \left(\frac{m\pi x_{2}}{a}\right)\sin \left( \omega t\right),  \label{psi2} \\
W_{(A)}^{nm} &=& A_{3}\sin \left( \frac{n\pi x_{1}}{a}\right) \sin \left(\frac{m\pi x_{2}}{a}\right)\sin \left( \omega t\right),  \label{w} \\
\Omega_{3(A)}^{nm}&=& A_{4}\cos \left( \frac{n\pi x_{1}}{a}\right) \cos \left(\frac{m\pi x_{2}}{a}\right)\sin \left( \omega t\right),  \label{omega3} \\
\Omega _{1(A)}^{0,nm} &=& A_{5}\sin \left( \frac{n\pi x_{1}}{a}\right) \cos \left(\frac{m\pi x_{2}}{a}\right)\sin \left( \omega t\right),  \label{omega1} \\
\Omega_{2(A)}^{0,nm}&=& A_{6}\cos \left( \frac{n\pi x_{1}}{a}\right) \sin \left(\frac{m\pi x_{2}}{a}\right)\sin \left( \omega t\right),  \label{omega2} \\
W_{(A)}^{\ast ,nm} &=& A_{7}\sin \left( \frac{n\pi x_{1}}{a}\right) \sin \left(\frac{m\pi x_{2}}{a}\right)\sin \left( \omega t\right),  \label{ws} \\
\hat{\Omega}_{1(A)}^{0,nm}&=& A_{8}\sin \left( \frac{n\pi x_{1}}{a}\right) \cos \left(\frac{m\pi x_{2}}{a}\right)\sin \left( \omega t\right),  \label{omega01} \\
\hat{\Omega}_{2(A)}^{0,nm} &=& A_{9}\cos \left( \frac{n\pi x_{1}}{a}\right) \sin \left(\frac{m\pi x_{2}}{a}\right)\sin \left( \omega t\right), \\ \label{omega02}
\Psi_{1(B)}^{nm} &=& B_{1}\sin \left( \frac{n\pi x_{1}}{a}\right) \cos \left(\frac{m\pi x_{2}}{a}\right) \sin \left( \omega t\right),  \label{psi1} \\
\Psi_{2(B)}^{nm}&=& B_{2}\cos \left( \frac{n\pi x_{1}}{a}\right) \sin \left(\frac{m\pi x_{2}}{a}\right) \sin \left( \omega t\right),  \label{psi2} \\
W_{(B)}^{nm} &=& B_{3}\cos \left( \frac{n\pi x_{1}}{a}\right) \cos \left(\frac{m\pi x_{2}}{a}\right) \sin \left( \omega t\right),  \label{w} \\
\Omega_{3(B)}^{nm}&=& B_{4}\sin \left( \frac{n\pi x_{1}}{a}\right) \sin \left(\frac{m\pi x_{2}}{a}\right) \sin \left( \omega t\right),  \label{omega3} \\
\Omega_{1(B)}^{0,nm} &=& B_{5}\cos \left( \frac{n\pi x_{1}}{a}\right) \sin \left(\frac{m\pi x_{2}}{a}\right) \sin \left( \omega t\right),  \label{omega1} \\
\Omega_{2(B)}^{0,nm}&=& B_{6}\sin \left( \frac{n\pi x_{1}}{a}\right) \cos \left(\frac{m\pi x_{2}}{a}\right) \sin \left( \omega t\right),  \label{omega2} \\
W_{(B)}^{\ast ,nm} &=& B_{7}\cos \left( \frac{n\pi x_{1}}{a}\right) \cos \left(\frac{m\pi x_{2}}{a}\right) \sin \left( \omega t\right),  \label{ws} \\
\hat{\Omega}_{1(B)}^{0,nm}&=& B_{8}\cos \left( \frac{n\pi x_{1}}{a}\right) \sin \left(\frac{m\pi x_{2}}{a}\right) \sin \left( \omega t\right),  \label{omega01} \\
\hat{\Omega}_{2(B)}^{0,nm} &=& B_{9}\sin \left( \frac{n\pi x_{1}}{a}\right) \cos \left(\frac{m\pi x_{2}}{a}\right) \sin \left( \omega t\right),  \label{omega02}
\end{eqnarray*}
and $A_{i}$ and $B_{i}$ are constants.

We solve the eigenvalue problem by substituting the expressions (\ref{psi1}) -- (\ref{omega02}) into the system of equations (\ref{bending_system}). The obtained 9 sequencies of positive eigenfrequencies $\omega_{i}^{nm}$ are associated with the rotation of the middle plane ($\omega_{1}^{nm}$ and $\omega_{2}^{nm}$), flexural motion and its transverse variation ($\omega _{3}^{nm}$ and $\omega _{7}^{nm}$), micro rotatory inertia ($\omega_{4}^{nm}$, $\omega_{5}^{nm}$ and $\omega _{6}^{nm}$) and its transverse variation ($\omega_{8}^{nm}$ and $\omega _{9}^{nm}$) \cite{Steinberg2015}.

\begin{figure}[H]
\begin{center}
\includegraphics[width=3.0in]{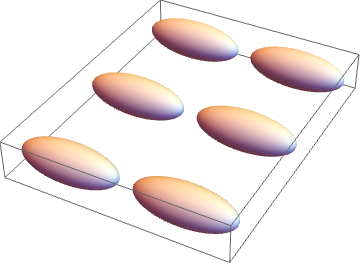} \vspace{.1in}
\end{center}
\caption{\textit{Horizontally stretched ellipsoid micro-elements: $J_{x}=0.002$, $J_{y}=0.001$, $J_{z}=0.0001$}}
\label{fig:ellipsoids}
\end{figure}

We perform all our numerical simulations for $a=3.0$m and $h=0.1$m. We consider different forms of micro elements: ball-shaped elements, horizontally and vertically stretched ellipsoids (see Figures \ref{fig:balls} and \ref{fig:ellipsoids}). For simplicity we will use the notation $\omega_{i}$ for the first elements $\omega_{i}^{11}$ of the sequences $\omega_{i}^{nm}$. The results of the computations are given in the Table \ref{tab:different}. 

The shape of the micro elements does not effect the natural macro frequencies $\omega_{1}$ and $\omega_{2}$ associated with the rotation of the middle
plane and $\omega_{3}$ and $\omega_{7}$ associated with the flexural motion and its transverse variation. 

\noindent
\begin{center}
\begin{table}[ht]
\caption{Eigenfrequencies $\omega _{i}^{11}$ (Hz) for different angles of rotation of horizontal ellipsoid micro-elements}
\label{tab:angle}%
\setkeys{Gin}{keepaspectratio} {
\begin{tabularx}{\textwidth}{@{}XXXXXXXXXX@{}}
\hline
Angle $\theta$ & $\omega_{1}$ & $\omega_{2}$ & $\omega_{3}$ & $\omega_{7}$ & $\omega_{4}$ & $\omega_{5}$ & $\omega_{8}$ & $\omega_{6}$ & $\omega_{9}$ \\ 
\hline
$0^{\circ}$ & 17.88 & 17.88 & 0.31 & 0.31  &	650.221 & 265.37 & 265.37 & 450.61 & 450.61 \\
$10^{\circ}$ & 17.88 & 17.88 & 0.31 & 0.31  &	650.221 & 255.59 & 279.40 & 429.89 & 469.93 \\
$20^{\circ}$ & 17.88 & 17.88 & 0.31 & 0.31  &	650.221 & 247.75 & 295.33 & 406.70 & 484.79 \\
$30^{\circ}$ & 17.88 & 17.88 & 0.31 & 0.31  &	650.221 & 242.57 & 313.65 & 382.94 & 495.14 \\ 
$40^{\circ}$ & 17.88 & 17.88 & 0.31 & 0.31  &	650.221 & 239.99 & 333.10 & 360.57 & 500.46 \\
$45^{\circ}$ & 17.88 & 17.88 & 0.31 & 0.31  &	650.221 & 239.68 & 338.95 & 354.35 & 501.13 \\
$50^{\circ}$ & 17.88 & 17.88 & 0.31 & 0.31  &	650.221 & 239.99 & 333.10 & 360.57 & 500.46 \\
$60^{\circ}$ & 17.88 & 17.88 & 0.31 & 0.31  &	650.221 & 242.57 & 313.65 & 382.94 & 495.14 \\
$70^{\circ}$ & 17.88 & 17.88 & 0.31 & 0.31  &	650.221 & 247.75 & 295.33 & 406.70 & 484.79 \\
$80^{\circ}$ & 17.88 & 17.88 & 0.31 & 0.31  &	650.221 & 255.59 & 279.40 & 429.89 & 469.93 \\
$90^{\circ}$ & 17.88 & 17.88 & 0.31 & 0.31  &	650.221 & 265.37 & 265.37 & 450.61 & 450.61 \\
\hline
\end{tabularx}
}
\end{table}
\end{center}

The ellipsoid elements have higher micro frequencies associated with the micro rotatory inertia ($\omega_{4}$, $\omega_{5}$ and $\omega_{6}$) and its transverse variation ($\omega_{8}$ and $\omega_{9}$), than the ball-shaped elements.

Let $J_{x}$, $J_{y}$ and $J_{z}$ be the principal moments of inertia of the
micro elements corresponding to the principal axes of their rotation. We
assume that the quantities $J_{x}$, $J_{y}$ and $J_{z}$ are constant
throughout the plate $B_{0}$. If the micro elements are rotated around the $%
z $-axis by the angle $\theta$ the rotatory inertia tensor $\mathbf{J}$ can
be expressed as 
\begin{equation}
\mathbf{J} = \left(%
\begin{matrix}
J_{x}\cos^{2}{\theta}+J_{y}\sin^{2}{\theta} & \left(J_{x}-J_{y}\right) \sin{%
2\theta} & 0 \\ 
\left(J_{x}-J_{y}\right) \sin{2\theta} & J_{x}\sin^{2}{\theta}+J_{y}\cos^{2}{%
\theta} & 0 \\ 
0 & 0 & J_{z}%
\end{matrix}
\right)
\end{equation}

The eigenfrequencies for different angles of rotation of the micro-elements are given in the Table \ref{tab:angle} and the Figure \ref{fig:microfreq}. The rotatory inertia principle
moments used are $J_{x}=0.002$, $J_{y}=0.001$, $J_{z}=0.0001$, which represent a horizontally stretched ellipsoid micro element (Figure \ref{fig:ellipsoids}). The case when the micro elements are not aligned with the edges of the plate the model predicts some additional natural frequencies related with the microstructure of the material.

\subsection{Three-dimensional Cosserat Body \\Vibration}

Let us consider the plate $B_{0}$ being a rectangular cuboid $[0,a]\times
\lbrack 0,a]\times \left[ -\frac{h}{2},\frac{h}{2} \right]$. Let the sets $T$
and $B$ be the top and the bottom surfaces contained in the planes $x_{3}=%
\frac{h}{2}$ and $x_{3}=-\frac{h}{2}$ respectively, and the curve $\Gamma =
\Gamma_{1} \cup \Gamma_{2}$ be the lateral part of the boundary: 
\begin{eqnarray*}
\Gamma_{1} = \left\{ \left( x_{1},x_{2},x_{3}\right) :x_{1} \in
\left\{0,a\right\} , x_{2}\in \left[ 0,a\right], x_{3}\in \left[ -\frac{h}{2}%
,\frac{h}{2} \right] \right\},  \notag \\
\Gamma_{2} = \left\{ \left( x_{1},x_{2},x_{3}\right) :x_{1}\in \left[ 0,a%
\right] , x_{2} \in \left\{0,a\right\}, x_{3}\in \left[ -\frac{h}{2},\frac{h%
}{2} \right] \right\},  \notag
\end{eqnarray*}

We solve the three-dimensional Cosserat equilibrium equations (\ref%
{equilibrium_equations 0}) -- (\ref{equilibrium_equations}) accompanied by
the constitutive equations (\ref{Hooke's_law 1}) -- (\ref{Hooke's_law 1A})
and strain-displacement and torsion-rotation relations (\ref{kinematic
formulas}) complemented by the following boundary conditions: 
\begin{eqnarray}
\Gamma _{1} &:& u_{2}=0,\text{ }u_{3}=0,\text{ }\varphi _{1}=0,
\label{BC3D_1} \\
\Gamma _{1} &:& \sigma_{11}=0,\text{ }\mu _{12}=0,\text{ }\mu _{13}=0;
\label{BC3D_2} \\
\Gamma _{2} &:& u_{1}=0,\text{ }u_{3}=0,\text{ }\varphi _{2}=0,
\label{BC3D_3} \\
\Gamma _{2} &:& \sigma_{22}=0,\text{ }\mu _{21}=0,\text{ }\mu _{23}=0;
\label{BC3D_4} \\
T &:&\sigma _{33}= p\left( x_{1}, x_{2} \right),\text{ }\mu _{33}=0;
\label{BC3D_5} \\
B &:&\sigma _{33}=0,\text{ }\mu_{33}=0.  \label{BC3D_6}
\end{eqnarray}

\noindent 
\begin{figure}[H]
\begin{center}
\includegraphics[width=3.0in]{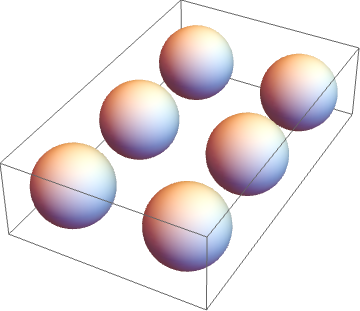} \vspace{.1in}
\end{center}
\caption{\textit{Ball-shaped micro-elements: $J_{x}=0.001$, $J_{y}=0.001$, $J_{z}=0.001$ }}
\label{fig:balls}
\end{figure}

\hspace{0.5in}

\noindent 

\noindent
\begin{center}
\begin{table}[H]
\caption{Comparison of the eigenfrequencies $\omega_{i}$ (Hz) with the exact values of the 3D Cosserat elasticity}
\label{tab:3d}
\setkeys{Gin}{keepaspectratio} {
\begin{tabularx}{\textwidth}{@{}p{.25\textwidth}XXXXX@{}}
\hline
 & $\omega_{1}$, $\omega_{2}$, & $\omega_{3}$, $\omega_{7}$ & $\omega_{4}$ & $\omega_{5}$, $\omega_{8}$ & $\omega_{6}$, $\omega_{9}$  \\ 
\hline
Plate theory & 0.310 & 17.881 & 501.13 & 205.62 & 338.95 \\
3D Cosserat elasticity & 0.309 & 17.763 & 530.82 & 211.98 & 317.87  \\
\hline
\end{tabularx}
}
\end{table}
\end{center}

Here the initial distribution of the pressure is given as 
\begin{equation*}
p =\sin \left( \frac{\pi x_{1}}{a}\right) \sin\left( \frac{\pi x_{2}}{a}\right) \sin{\omega t}
\end{equation*}
and the rotatory inertia tensor $\mathbf{J}$ is assumed to have a diagonal
form 
\begin{equation}
\mathbf{J} = \left[
\begin{matrix}
J_{x} & 0 & 0 \\ 
0 & J_{y} & 0 \\ 
0 & 0 & J_{z}
\end{matrix}
\right].
\end{equation}

Using the method of separation of variables and taking into account the boundary conditions (\ref{BC3D_1}) -- (\ref{BC3D_4}), we express the kinematic variables in the form: 
\begin{eqnarray}
&&u_{1} =\cos \left( \frac{\pi x_{1}}{a}\right) \sin \left( \frac{\pi x_{2}}{%
a}\right) z_{1}\left( x_{3}\right) \sin{\omega t} ,
\label{solutionexpressions3Daa} \\
&& u_{2}=\sin \left( \frac{\pi x_{1}}{a}\right) \cos \left( \frac{\pi x_{2}}{%
a}\right) z_{2}\left( x_{3}\right) \sin{\omega t} ,
\label{solutionexpressions3Da} \\
&& u_{3} =\sin \left( \frac{\pi x_{1}}{a}\right) \sin \left( \frac{\pi x_{2}%
}{a}\right) z_{3}\left( x_{3}\right) \sin{\omega t},
\label{solutionexpressions3Dbc} \\
&& \phi _{1}=\sin \left( \frac{\pi x_{1}}{a}\right) \cos \left( \frac{\pi
x_{2}}{a} \right) z_{4}\left( x_{3}\right) \sin{\omega t} ,
\label{solutionexpressions3Dbb} \\
&& \phi _{2} = \cos \left( \frac{\pi x_{1}}{a}\right) \sin \left( \frac{\pi
x_{2}}{a}\right) z_{5}\left( x_{3}\right) \sin{\omega t},
\label{solutionexpressions3Dg} \\
&& \phi _{3}=\cos \left( \frac{\pi x_{1}}{a}\right) \cos \left( \frac{\pi
x_{2}}{a}\right) z_{6}\left( x_{3}\right) \sin{\omega t},
\label{solutionexpressions3Db}
\end{eqnarray}
where the functions $z_{i}\left( x_{3}\right)$ represent the transverse
variations of the kinematic variables.

If we substitute the expressions (\ref{solutionexpressions3Daa}) -- (\ref{solutionexpressions3Db}) into (\ref{Hooke's_law 1}) -- (\ref{Hooke's_law 1A}) and then into (\ref{equilibrium_equations 0}) -- (\ref{equilibrium_equations}) we will obtain the following eigenvalue problem 
\begin{equation}
\mathbf{B} z = \omega^{2} \mathbf{A} z \label{3d_system}
\end{equation}
where
\begin{equation*}
\mathbf{B} = L_{0} \mathbf{C_{0}} + L_{1} \mathbf{C_{1}} + L_{2} \mathbf{C_{2}}
\end{equation*}
and
\begin{equation}
\mathbf{C_{0}} = 
\begingroup 
\setlength\arraycolsep{2pt}
\begin{pmatrix}
b_{2} & b_{3} & 0 & 0 & 0 & -b_{6} \\ 
b_{3} & b_{2} & 0 & 0 & 0 & b_{6} \\ 
0 & 0 & 0 & b_{6} & -b_{6} & 0 \\ 
0 & 0 & b_{6} & b_{10} & b_{11} & 0 \\ 
0 & 0 & -b_{6} & b_{11} & b_{10} & 0 \\ 
-b_{6} & b_{6} & 0 & 0 & 0 & b_{14}
\end{pmatrix}
\endgroup
\end{equation}

\begin{equation}
\mathbf{C_{1}} = 
\begingroup 
\setlength\arraycolsep{2pt}
\begin{pmatrix}
0 & 0 & b_{4} & 0 & -b_{5} & 0 \\
0 & 0 & b_{4} & b_{5} & 0 & 0 \\
-b_{4} & b_{4} & 0 & 0 & 0 & 0 \\
0 & b_{5} & 0 & 0 & 0 & b_{12} \\
b_{5} & 0 & 0 & 0 & 0 & b_{12} \\
0 & 0 & 0 & -b_{12} & b_{12} & 0
\end{pmatrix}
\endgroup
\end{equation}

\begin{equation}
\mathbf{C_{2}} = 
\begingroup 
\setlength\arraycolsep{2pt}
\begin{pmatrix}
b_{1} & 0 & 0 & 0 & 0 & 0 \\ 
0 & b_{1} & 0 & 0 & 0 & 0 \\ 
0 & 0 & b_{7} & 0 & 0 & 0 \\ 
0 & 0 & 0 & b_{9} & 0 & 0 \\ 
0 & 0 & 0 & 0 & b_{9} & 0 \\ 
0 & 0 & 0 & 0 & 0 & b_{13}
\end{pmatrix}
\endgroup
\end{equation}

\begin{equation}
\mathbf{A} =
\begingroup 
\setlength\arraycolsep{2pt}
\begin{pmatrix}
-a^{2}\rho & 0 & 0 & 0 & 0 & 0 \\ 
0 & -a^{2}\rho & 0 & 0 & 0 & 0 \\ 
0 & 0 & -a^{2}\rho & 0 & 0 & 0 \\ 
0 & 0 & 0 & -a^{2}J_{x} & 0 & 0 \\ 
0 & 0 & 0 & 0 & -a^{2}J_{y} & 0 \\ 
0 & 0 & 0 & 0 & 0 & -a^{2}J_{z} \\ 
\end{pmatrix}
\endgroup
\end{equation}
\begin{equation}
z = 
\begingroup 
\setlength\arraycolsep{2pt}
\begin{pmatrix}
z_{1}, & z_{2}, & z_{3}, & z_{4}, & z_{5}, & z_{6}
\end{pmatrix}
\endgroup^{T}
\end{equation}
\noindent 
\begin{figure}[H]
\begin{center}
\includegraphics[width=3.5in]{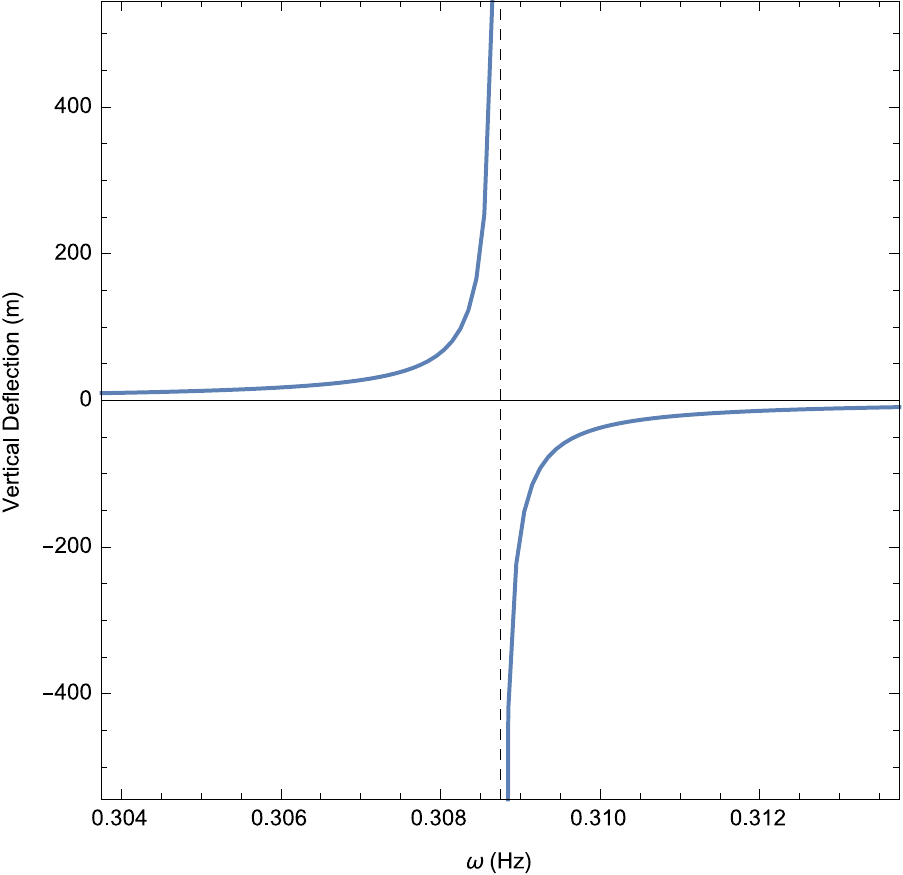}
\end{center}
\caption{\textit{Large amplitude vibrations of the Cosserat body forced to vibrate at its natural frequency $\omega_{1}$.}}
\label{fig:0309}
\end{figure}

and the differential operators $L_{i}$ are defined as
\begin{equation*}
L_{0} = I, \hspace{0.2in} L_{1} = \frac{d}{dx_{3}}, \hspace{0.2in} L_{2} =\frac{d^{2}}{d x^{2}_{3}}
\end{equation*}
and the coefficients $b_{i}$ are defined as
\begin{align*}
b_{1} & = a^{2}\left( \mu +\alpha \right), & b_{2} & = -\pi ^{2}(\alpha+\lambda+3\mu), \\ b_{3} & = -\pi ^{2}(\lambda +\mu -\alpha), & b_{4} & = a\pi (\lambda +\mu -\alpha ), \\
b_{5} & = 2a^{2}\alpha, & b_{6} & = 2a\pi \alpha, \\
b_{7} & = a^{2}(2\mu +\lambda), &b_{8} & = -2\pi^{2}(\alpha+\mu), \\ 
b_{9} & = a^{2}(\gamma +\epsilon), &
b_{10} & = -\pi^{2}(\beta+\epsilon+3\gamma), \\
b_{11} & = -\pi^{2}(\beta+\gamma-\epsilon), &
b_{12} & = -a\pi(\beta+\gamma-\epsilon), \\
b_{13} & = a^{2}(\beta+2\gamma), & b_{14} & = -2\pi^{2}(\gamma+\epsilon) -4a^2\alpha
\end{align*}
The system of differential equations (\ref{3d_system}) is complemented by the following boundary conditions 
\begin{eqnarray}
x_{3}=\frac{h}{2}: && \mathbf{D} z = D_{0}, \\
x_{3}=-\frac{h}{2}: && \mathbf{D} z = 0,
\end{eqnarray}
where
\begin{equation}
\mathbf{D} = 
\begingroup 
\setlength\arraycolsep{2pt}
\begin{pmatrix}
d_{1}L_{1} & 0 & d_{2}L_{0} & 0 & -d_{3}L_{0} & 0 \\ 
0 & d_{1}L_{1} & d_{2}L_{0} & d_{3}L_{0} & 0 & 0 \\ 
d_{4}L_{0} & d_{4}L_{0} & d_{5}L_{1} & 0 & 0 & 0 \\ 
0 & 0 & 0 & d_{6}L_{1} & 0 & d_{7}L_{0} \\ 
0 & 0 & 0 & 0 & d_{6}L_{1} & d_{7}L_{0} \\ 
0 & 0 & 0 & d_{8}L_{0} & d_{8}L_{0} & d_{9}L_{1} \\ 
\end{pmatrix}
\endgroup
\end{equation}
\begin{equation}
D_{0}=
\begingroup 
\setlength\arraycolsep{2pt}
\begin{pmatrix}
0, & 0, & a, & 0, & 0, & 0 
\end{pmatrix}
\endgroup^{T},
\end{equation}
and the coefficients $d_{i}$ are defined as
\begin{align*}
d_{1} & = a(\mu +\alpha), & d_{2} & = -\pi (\mu -\alpha), \\ 
d_{3} & = 2a\alpha, & d_{4} & = a(\lambda +2\mu ), \\
d_{5} & = -\pi \lambda, & d_{6} & = a(\gamma +\epsilon ), \\
d_{7} & = a(\gamma -\epsilon ), & d_{8} & = \pi\beta, \\ 
d_{9} & = a(\beta +2\gamma ).
\end{align*}

We solve this eigenvalue problem by forcing the Cosserat body to vibrate at a given frequency $\omega$. When the frequency $\omega$ coincides with the natural frequency of the plate the resonance will occur and the large amplitude vibrations can be observed (Figure \ref{fig:0309}).

The comparison of the eigenfrequencies of the Cosserat plate with the eigenfrequencies of the three-dimensional Cosserat elasticity is given in the Table \ref{tab:3d}. The rotatory inertia principle moments used are $J_{x}=0.001$, $J_{y}=0.001$, $J_{z}=0.001$, which represent a ball-shaped micro element (Figure \ref{fig:balls}). The relative error of the natural macro frequencies associated with the rotation of the middle plane and the flexural motion is less than 1\%.

\section{Conclusion}

In this paper we presenedt the validation of our mathematical model for the dynamics of Cosserat elastic plates. The validation of the model was based on the comparison with the exact solution of the 3-dimensional Cosserat elastodynamics. The computations of eigenfrequencies show the high agreement with the exact values. This allowed us to detect the splitting of the frequencies of vibrations (micro vibration) depending on the orientation of micro elements. We showed that this approach is a powerful tool for distinguishing between the frequencies of the micro and macro vibrations of the plate.

\end{document}